\documentstyle[12pt]{article}
\catcode`\@=11
\@addtoreset{equation}{section}

\catcode`\@=12
\newtheorem{Theorem}{Theorem}[section]
\newtheorem{Definition}{Definition}[section]
\newtheorem{Proposition}{Proposition}[section]
\newtheorem{Lemma}{Lemma}[section]
\newtheorem{Corollary}{Corollary}[section]
\newtheorem{Note}{Note}[section]
\title{$J-$holomorphic disks and Lagrangian Squeezing In 
$M\times C$\thanks{Project 19871044 Supported by NSF}}

\author{Renyi Ma\\
Department of Mathmatics \\
Tsinghua University \\
Beijing, 100084\\
People's Republic of China}

\date { }

\begin{document}
\textwidth=125mm
\textheight=185mm
\parindent=8mm
\frenchspacing
\maketitle

\begin{abstract}
In this article, we define an invariant for Lagrangian 
submanifold and  prove that if the Lagrangian submanifold 
contained in the ball of radius $r$, then the invariant is less than $4\pi r^2$. This 
modifies Gromov's Lagrangian embedding theorem. 
\end{abstract}
\noindent{\bf Keywords} Symplectic geometry, J-holomorphic curves, 
Chord.

\noindent{\bf 2000MR Subject Classification} 32Q65,53D35,53D12

\section{Introduction and results}

Let $(M,\omega )$ be a symplectic manifold.
Let $J$ be the almost complex structure tamed by $\omega $, i.e., 
$\omega (v, Jv)>0$ for $v\in TM$. Let 
${\cal {J}}$ the space of all tame almost complex structures. 
 
\begin{Definition}
Let 
$$s(M,\omega ,J)=\inf \{ \int _{S^2}f^*\omega >0 |f:S^2\to M \ is \ J-holomorphic \} $$
\end{Definition}

\begin{Definition}
Let 
$$s(M,\omega )=\sup _{J\in {\cal {J}}}l(M,\omega ,J)$$
\end{Definition}

Let $W$ be a Lagrangian submanifold in $M$, i.e., 
$\omega |W=0$. 

\begin{Definition}
Let 
$$l(M,W, \omega )=\inf \{ |\int _{D^2}f^*\omega |>0 |f:(D^2,\partial D^2)\to (M,W) \} $$
\end{Definition}

\begin{Theorem}
Let $(M,\omega )$ be a closed  
compact symplectic manifold or a manifold convex at infinity and 
$M\times C$ be a symplectic manifold with 
symplectic form $\omega \oplus \sigma $, 
here $(C,\sigma )$ standard 
symplectic plane. Let 
$4\pi r_0^2<s(M,\omega )$ and 
$B_{r_0}(0)\subset C$ the closed disk with radius 
$r_0$. If $W$ is a close Lagrangian manifold in $M\times B_{r_0}(0)$, then 
$$l(M,W, \omega )<4\pi r_0^2$$
\end{Theorem}
This can be considered as an Lagrangian version of Gromov's symplectic squeezing.

\begin{Corollary}
(Gromov\cite{gro})Let $(V',\omega ')$ be an exact symplectic manifold with 
restricted contact boundary and $\omega '=d\alpha ' $. Let 
$V'\times C$ be a symplectic manifold with 
symplectic form $\omega '\oplus \sigma =d\alpha 
=d(\alpha '\oplus \alpha _0$, 
here $(C,\sigma )$ standard 
symplectic plane. 
If $W$ is a close exact Lagrangian submanifold, then 
$l(V'\times C,W, \omega )==\infty $, i.e., there does not exist any 
close exact Lagrangian submanifold in $V'\times C$.
\end{Corollary}

\begin{Corollary}
Let $L^n$ be a close Lagrangian in $R^{2n}$ and $L(R^{2n},L^n,\omega )=4\pi r_0^2>0$, then 
$L^n$ can not be embedded in $B_{r_0}(0)$ as a Lagrangian submanifold. 
\end{Corollary}

{\bf Sketch of proofs:}We will work in the framework 
proposed by Gromov in \cite{gro}. In Section 2, we study the 
linear Cauchy-Riemann operator and sketch some basic 
properties. In section 3, we study 
the space ${\cal D}(V,W)$
of disks in manifold $V$ with boundary 
in Lagrangian submanifold $W$. 
In Section 3, we construct a Fredholm
section of tangent bundle of  ${\cal  D }(V,W)$.  
In Section 4, we construct an anti-holomorphic section 
and prove that it is not proper. 
In section 5, 
we
use Gromov's nonlinear Fredholm alternative to 
prove Theorem 1.1.

\section{Linear Fredholm theory}

For $100<k<\infty $ consider the Hilbert space
$V_k$ consisting of all maps $u\in H^{k,2}(D, C^n)$,
such that $u(z)\in R^n\subset C^n$ for almost all $z\in
\partial D$. $L_{k-1}$ denotes the usual Hilbert $L_{k-1}-$space
$H_{k-1}(D, C^n)$. We define an operator
$\bar \partial :V_p\mapsto L_p$ by
\begin{equation}
\bar \partial u=u_s+iu_t
\end{equation}
where the coordinates on $D$ are
$(s,t)=s+it$, $D=\{ z||z|\leq 1\}  $.
The following result is well known, see\cite{wen}.
\begin{Proposition}
$\bar \partial :V_p\mapsto L_p $ is a surjective real
linear Fredholm operator of index $n$. The kernel
consists of the constant real valued maps.
\end{Proposition}
Let $(C^n, \sigma =-Im(\cdot ,\cdot ))$ be the standard
symplectic space. We consider a real $n-$dimensional plane
$R^n\subset C^n$. It is called Lagrangian if the
skew-scalar product of any two vectors of $R^n$ equals zero.
For example, the plane $p=0$ and $q=0$ are Lagrangian
subspaces. The manifold of all (nonoriented) Lagrangian subspaces of
$R^{2n}$ is called the Lagrangian-Grassmanian $\Lambda (n)$.
One can prove that the fundamental group of
$\Lambda (n)$ is free cyclic, i.e.
$\pi _1(\Lambda (n))=Z$. Next assume
$(\Gamma (z))_{z\in \partial D}$ is a smooth map
associating to a point $z\in \partial D$ a Lagrangian
subspace $\Gamma (z)$ of $C^n$, i.e.
$(\Gamma (z))_{z\in \partial D}$ defines a smooth curve
$\alpha $ in the Lagrangian-Grassmanian manifold $\Lambda (n)$.
Since $\pi _1(\Lambda (n))=Z$, one have
$[\alpha ]=ke$, we call integer $k$ the Maslov index
of curve $\alpha $ and denote it by $m(\Gamma )$, see(\cite{ag}).

Now let $z:S^1\mapsto R^n\subset C^n$ be a smooth curve. 
Then
it defines a constant loop $\alpha $ in Lagrangian-Grassmanian
manifold $\Lambda (n)$. This loop defines
the Maslov index $m(\alpha )$ of the map
$z$ which is easily seen to be zero.

  Now Let $(V,\omega )$ be a symplectic manifold and 
$W\subset V$ a closed Lagrangian submanifold. Let 
$u:D^2\to V$ be a smooth map homotopic to constant map 
with boundary 
$\partial D\subset W$.
Then $u^*TV$ is a symplectic vector bundle and 
$(u|_{\partial D})^*TW$ be a Lagrangian subbundle in 
$u^*TV$. Since $u$ is contractible, we can take 
a trivialization of $u^*TV$ as 
$$\Phi (u^*TV)=D\times C^n$$ 
and 
$$\Phi (u|_{\partial D})^*TW)\subset S^1\times C^n$$
Let 
$$\pi _2: D\times C^n\to C^n$$
then 
$$\bar u: z\in S^1\to \pi _2\Phi (u|_{\partial D})^*TW(z)\in \Lambda (n).$$
Write $\bar u=u|_{\partial D}$.
\begin{Lemma}
Let $u: (D^2,\partial D^2) \rightarrow (V,W)$ be a $C^k-$map $(k\geq 1)$ 
as above. Then,
$$m(u|_{\partial D})=0$$
\end{Lemma}
Proof.  Since $u$ is contractible in $V$ relative to $W$, we have 
a homotopy $\Phi _s$ of trivializations
such that 
$$\Phi _s(u^*TV)=D\times C^n$$ 
and 
$$\Phi _s(u|_{\partial D})^*TW)\subset S^1\times C^n$$
Moreover 
$$\Phi _0(u|_{\partial D})^*TW=S^1\times R^n$$
So, the homotopy induces 
a homotopy $\bar h$ in Lagrangian-Grassmanian
manifold. Note that $m(\bar h(0, \cdot ))=0$.
By the homotopy invariance of Maslov index,
we know that $m(u|_{\partial D})=0$.

\vskip 5pt 

   Consider the partial differential equation
\begin{eqnarray}
\bar \partial u+A(z)u=0  \ on \ D  \\
u(z)\in \Gamma (z) R^n\ for \ z\in \partial D \\
\Gamma (z)\in GL(2n,R)\cap Sp(2n)\\
m(\Gamma )=0 \ \ \ \ \ \ \ \ 
\end{eqnarray}

For $100<k<\infty $ consider the Banach space $\bar V_k $
consisting of all maps $u\in H^{k,2}(D, C^n)$ such 
that  $u(z)\in \Gamma (z)$ for almost all $z\in
\partial D$. Let $L_{k-1}$ the usual $L_{k-1}-$space $H_{k-1}(D,C^n)$ and 

$$L_{k-1}(S^1)=\{ u\in H^{k-1}(S^1)|u(z)\in \Gamma (z) R^n\ for \ z\in 
\partial D\}$$ 
We define an operator $P$:
$\bar V_{k}\rightarrow L_{k-1}\times L_{k-1}(S^1)$ by
\begin{equation}
P(u)=(\bar \partial u+Au,u|_{\partial D})
\end{equation}
where $D$ as in (2.1).
\begin{Proposition}
$\bar \partial : \bar V_p \rightarrow L_p$
is a real linear Fredholm operator of index n.
\end{Proposition}

\section{Nonlinear Fredholm Theory}

Let $(M,\omega )$ be a closed symplectic 
manifold with 
almost complex structure $J_M$ tamed by $\omega $, i.e.,  
$\omega (\cdot ,J_M\cdot )$ determines 
a Riemannian metric $g_M$. Let $(R^2,\omega _0)$ 
be the standard symplectic vector space  
with the adapted metric $g_0 $ and $J_0=i$.
Let $W$ is a close Lagrangian submanifold in $M\times C$.

\vskip 5pt

   Let 
$${\cal D}^k(V,W,p)=\{ u \in H^k(D,V)|
u(x)\in W \ a.e \ for \ x\in \partial D \ and \ u(1)=p\}$$
for $k\geq 100$.
\begin{Lemma}
Let $W$ be a closed Lagrangian submanifold in 
$V$. Then, 
$${\cal D}^k(V,W,p)=\{ u \in H^k(D,V)|
u(x)\in W \ a.e \ for \ x\in \partial D \ and \ u(1)=p\}$$
is a pseudo-Hilbert manifold with the tangent bundle
\begin{equation}
T{\cal D}^k(V,W,p)=\bigcup _{u\in {\cal {D}}^k(V,W,p)}
\Lambda ^{k-1}(u^*TV,u|_{\partial D}^*TW,p)
\end{equation}
here 
$$\Lambda ^{k-1}(u^*TV,u|_{\partial D}^*TW,p)=$$ 
$$\{ H^{k-1}-sections \ of \ (u^*(TV),(u|_{\partial D})^*TL)\ 
which \ vanishes  \ at \ 1\} $$
\end{Lemma}
\begin{Note}
Since $W$ is not regular we know that ${\cal D}^k(V,W,p)$ is 
in general complete, however it is enough for our purpose.
\end{Note}
Proof: See \cite{al,kl}. 

\vskip 3pt

   Now we consider  a section
from ${\cal D}^k(V,W,p)$ to $T{\cal D}^k(V,W, p)$ follows as in 
\cite{al,gro}, i.e., 
let $\bar \partial :{\cal D}^k(V,W,p)\rightarrow T{\cal D}^k(V,W,p)$
be the Cauchy-Riemmann section 
\begin{equation}
\bar \partial u={{\partial u}\over {\partial s}}
+J{{\partial u}\over {\partial t}}  \label{eq:CR}
\end{equation}
for $u\in {\cal D}^k(V,W,p)$.

\begin{Theorem}
The Cauchy-Riemann section $\bar \partial $ defined in (\ref{eq:CR})
is a Fredholm section of Index zero.
\end{Theorem}
Proof. According to the definition of the Fredholm section, 
we need to prove that
$u\in {\cal D}^k(V,W,p)$, the linearization
$D\bar \partial (u)$ of $\bar \partial $ at $u$ is a linear Fredholm
operator.
Note that
\begin{equation}
D\bar \partial (u)=D{\bar \partial _{[u]}}
\end{equation}
where
\begin{equation}
(D\bar \partial _{[u]})v=\frac{\partial v}{\partial s}
+J\frac{\partial v}{\partial t}+A(u)v
\end{equation}
with 
$$v|_{\partial D}\in (u|_{\partial D})^*TW$$
here $A(u)$ is $2n\times 2n$
matrix induced by the torsion of
almost complex structure, see \cite{al,gro} for the computation.

   Observe that the linearization $D\bar \partial (u)$ of 
$\bar \partial $ at $u$ is equivalent to the following Lagrangian 
boundary value problem
\begin{eqnarray}
&&{{\partial v}\over {\partial s}}+J{{\partial v}\over {\partial t}}
+A(u)v=f, \ v\in \Lambda ^k(u^*TV)\cr 
&&v(t)\in T_{u(t)}W, \ \ t\in {\partial D}  \label{eq:Lin}
\end{eqnarray}
One 
can check that (\ref{eq:Lin}) 
defines a linear Fredholm operator. In fact, 
by proposition 2.2 and Lemma 2.1, since the operator $A(u)$ is a compact, 
we know that the operator $\bar \partial $ is a nonlinear Fredholm operator 
of the index zero.

\begin{Definition}
Let $X$ be a Banach manifold and $P:Y\to X$ the Banach 
vector bundle.
A Fredholm section $F:X\rightarrow Y$ is
proper if $F^{-1}(0)$ is a compact set and is called 
generic if $F$ intersects the zero section transversally, see \cite{al,gro}.
\end{Definition}
\begin{Definition}
$deg(F,y)=\sharp \{ F^{-1}(0)\} mod2$ is called the Fredholm
degree of a Fredholm section (see\cite{al,gro}).
\end{Definition}
\begin{Theorem}
Assum that the Fredholm section
$F=\bar \partial : {\cal D}^k(V,W,p)\rightarrow T^({\cal D}^k(V,W,p)$
constructed in (\ref{eq:CR}) is proper. Then,
$$deg(F,0)=1$$
\end{Theorem}
Proof: We assume that $u:D\mapsto V$ be a $J-$holomorphic disk
with boundary $u(\partial D)\subset W$ and 
by the assumption that $u$ is homotopic to the 
constant map $u_0(D)=p$. Since almost complex
structure ${J}$ tamed by  the symplectic form $\omega $,
by stokes formula, we conclude $u: D\rightarrow
V$ is a constant map. Because $u(1)=p$, We know that
$F^{-1}(0)={p}$.
Next we show that the linearizatioon $DF(p)$ of $F$ at $p$ is
an isomorphism from $T_p{\cal D}(V,W,p)$ to $E$.
This is equivalent to solve the equations
\begin{eqnarray}
{\frac {\partial v}{\partial s}}+J{\frac {\partial v}{\partial t}}
+Av=f\\
v|_{\partial D}\subset T_pW
\end{eqnarray}
here $J=J(p)=i$ and $A$ a constant matrix. By Lemma 2.1, we know that $DF(p)$ is an isomorphism.
Therefore $deg(F,0)=1$.

\section{Non-properness of a Fredholm section}

In this section we shall construct a non-proper Fredholm section 
$F_1:{\cal D}\rightarrow E$ by perturbing 
the Cauchy-Riemann section as in 
\cite{al,gro}.

\subsection{Anti-holomorphic section}

  Let $(V',\omega ')=(M,\omega _M)$ and  
$(V,\omega )=(V'\times C, \omega '\oplus \omega _0)$, 
for convenience we assume $n=1$, and 
$W$ as in section3 and $J=J'\oplus i$, $g=g'\oplus g_0$, 
$g_0$ the standard metric on $C$. 

   Now let $c\in C$ be a non-zero vector. We consider the 
equations
\begin{eqnarray}
v=(v',f):D\to V'\times C \nonumber \\
\bar \partial _{J'}v'=0,\bar \partial f=c\ \ \ \nonumber \\
v|_{\partial D}:\partial D\to W\ \ \ 
\end{eqnarray}
here $v$ homotopic to constant map 
$\{ p\}$ relative to $W$. 
Note that $W\subset V\times B_{r_0}(0)$ for $4\pi r_0^2<s(M,\omega )$. 
\begin{Lemma}
Let $v$ be the solutions of (4.1), then one has 
the following estimates
\begin{eqnarray}
E({v})=
\int _D(g'({{\partial {v'}}\over {\partial x}},
{J'}{{\partial {v'}}\over {\partial x}})
+g'({{\partial {v'}}\over {\partial y}},
{J'}{{\partial {v'}}\over {\partial y}}) \nonumber \\
+g_0({{\partial {f}}\over {\partial x}},
{i}{{\partial {f}}\over {\partial x}})
+g_0({{\partial {f}}\over {\partial y}},
{i}{{\partial {f}}\over {\partial y}}))d\sigma 
\leq 4\pi r_0^2 <s(M,\omega ). 
\end{eqnarray}
\end{Lemma}
Proof: Since $v(z)=(v'(z),f(z))$ satisfy (4.1)
and $v(z)=(v'(z),f(z))\in V'\times C$ 
is homotopic to constant map $v_0:D\to \{ p\}\subset W$ 
in $(V,W)$, by the Stokes formula
\begin{equation}
\int _{D}v^*(\omega '\oplus \omega _0)=0
\end{equation}
Note that the metric $g$ is adapted to the symplectic form 
$\omega $ and $J$, i.e., 
\begin{equation}
g=\omega  (\cdot ,J\cdot )
\end{equation}
By the simple algebraic computation, we have 
\begin{equation}
\int _{D}{v}^*\omega  ={{1}\over {4}}
\int _{D^2}(|\partial v|^2 
-|\bar {\partial }v|^2)=0
\end{equation}
and 
\begin{equation}
|\nabla v|={{1}\over {2}}(
|\partial v|^2 +|\bar \partial v|^2 
\end{equation}
Then 
\begin{eqnarray}
E(v)&=&\int _{D} |\nabla v| \nonumber \\ 
      &=&\int _{D}\{ {{1}\over {2}}(
|\partial v|^2+|\bar \partial v|^2)\} d\sigma \nonumber \\ 
&=&\pi |c|_{g_0}^2
\end{eqnarray}
By the equations (4.1), 
one get 
\begin{equation}
\bar \partial f=c \ on \ D
\end{equation}
We have 
\begin{equation}
f(z)={{1}\over {2}}c\bar z+h(z)
\end{equation}
here $h(z)$ is a holomorphic function on $D$. Note that  
$f(z)$ is smooth up to the boundary $\partial D$, then, by 
Cauchy integral formula
\begin{eqnarray}
\int _{\partial D}f(z)dz&=&{{1}\over {2}}c\int _{\partial D}
\bar {z}dz+\int _{\partial D}h(z)dz \cr
&=&\pi ic
\end{eqnarray}
So, we have 
\begin{equation}
|c|={{1}\over {\pi}}|\int _{\partial D^2}f(z)dz|
\end{equation}
Therefore, 
\begin{eqnarray}
E(v)&\leq &\pi |c|^2
\leq {{1}\over {\pi }}|\int _{\partial D}f(z)dz|^2      \cr
&\leq &{{1}\over {\pi }}|\int _{\partial D}|f(z)||dz||^2   \cr
&\leq &4\pi |diam(pr_2(W))^2 \cr
&\leq &4\pi r_0^2.
\end{eqnarray}
This finishes the proof of Lemma.

\begin{Proposition}
For $|c|\geq 3r_0$, then the 
equations (4.1)
has no solutions. 
\end{Proposition}
Proof. By (4.11), we have 
\begin{eqnarray} 
|c|&\leq &{{1}\over {\pi }}\int _{\partial D}|f(z)||dz|\cr
&\leq &{{1}\over {\pi }}\int _{\partial D}
diam(pr_2(W))||dz| \cr
&\leq &2r_0
\end{eqnarray}
It follows that $c=3r_0$ can not be obtained by 
any solutions.

\subsection{Modification of section $c$}

Note that the section $c$ is not a section of the 
Hilbert bundle in section 3 since $c$ is not 
tangent to the Lagrangian submanifold $W$, we must modify it as follows:

\vskip 3pt 

  Let $c$ as in section 4.1, we define 
\begin{eqnarray}
c_{\chi ,\delta }(z,v)=\left\{ \begin{array}{ll}
c \ \ \ &\mbox{if\  $|z|\leq 1-2\delta $,}\cr
0 \ \ \ &\mbox{otherwise}
\end{array}
\right. 
\end{eqnarray}
Then by using the cut off function $\varphi _h(z)$ and 
its convolution with section 
$c_{\chi ,\delta }$, we obtain a smooth section 
$c_\delta$ satisfying

\begin{eqnarray}
c_{\delta }(z,v)=\left\{ \begin{array}{ll}
c \ \ \ &\mbox{if\  $|z|\leq 1-3\delta $,}\cr
0 \ \ \ &\mbox{if\  $|z|\geq 1-\delta $.}
\end{array}
\right. 
\end{eqnarray}
for $h$ small enough, for the convolution theory see \cite{hor}.

   Now let $c\in C$ be a non-zero vector and 
$c_\delta $ the induced anti-holomorphic section. We consider the 
equations
\begin{eqnarray}
v=(v',f):D\to V'\times C \nonumber \\
\bar \partial _{J'}v'=0,\bar \partial f=c_\delta \ \ \ \nonumber \\
v|_{\partial D}:\partial D\to W\ \ \  \label{eq:4.16}
\end{eqnarray}
Note that $W\subset V\times B_{r_0}(0)$ for $4\pi r_0^2<s(M,\omega )$. 
Then by repeating the same argument as section 4.1., we obtain 
\begin{Lemma}
Let $v$ be the solutions of (\ref{eq:4.16}) and $\delta $ 
small enough, then one has 
the following estimates
\begin{eqnarray}
E({v})\leq 4\pi r_0^2 <s(M,\omega ). 
\end{eqnarray}
\end{Lemma}
and

\begin{Proposition}
For $|c|\geq 3r_0$, then the 
equations (\ref{eq:4.16})
has no solutions. 
\end{Proposition}

\subsection{Modification of $J\oplus i$}

Now let $V=M\times R^{2n}$ and $
J_1=J_M\oplus i$ on $V$ and $J_2=\bar J$ any almost complex structure 
on $M\times R^{2n}$. 
Now we consider the product $D^2\times V$ and 
almost complex structures on it.

\begin{eqnarray}
J_{\chi ,\delta }(z,v)=\left\{ \begin{array}{ll}
i\oplus J_M\oplus i \ \ \ &\mbox{if\  $|z|\leq 1-2\delta $,}\cr
i\oplus \bar J \ \ \ &\mbox{otherwise}
\end{array}
\right. 
\end{eqnarray}
Then by using the cut off function $\varphi _h(z)$ and 
its convolution with section 
$J_{\chi ,\delta }$, we obtain a smooth section 
$J_\delta$ satisfying

\begin{eqnarray}
J_{\delta }(z,v)=\left\{ \begin{array}{ll}
i\oplus J_M\oplus i \ \ \ &\mbox{if\  $|z|\leq 1-3\delta $,}\cr
i\oplus {\bar {J}} \ \ \ &\mbox{if\  $|z|\geq 1-\delta $.}
\end{array}
\right. 
\end{eqnarray}
for $h$ small enough, for the convolution theory see \cite{hor}.

Then as in section 4.2, one can also reformulation 
of the equations (\ref{eq:4.16}) and get similar 
estimates of Cauchy-Riemann equations, we leave it 
as exercises to reader.

\begin{Theorem}
The Fredholm sections $F_1=\bar \partial +c_\delta 
: {\cal  {D}}^k(V,W,p) \rightarrow T({\cal {D}}^k(V,W,p)$ is not proper 
for $|c|$ large enough.
\end{Theorem}
Proof. See \cite{al,gro}.

\section{$J-$holomorphic section}

In this section we 
show that the boundaries of the Cauchy-Riemann solutions 
of (\ref{eq:4.16}) remain  
in a finite part of the Lagrangian submanifold $W$ by 
using the monotone inequality.

\subsection{J-holomorphic section}

Recall that $W\subset \Sigma \subset 
V=\subset M\times R^{2n}$ as in section 3. 
The Riemann metric $g$ on $M\times R^{2n}$ 
induces a metric $g|W$. We assume that $n=1$ 
for convenience.

   Now let $c\in C$ be a non-zero vector and 
$c_\delta $ the induced anti-holomorphic section. We consider the 
nonlinear inhomogeneous equations (\ref{eq:4.16}) and 
transform it into $\bar J-$holomorphic map by 
considering its graph as in \cite{al,gro}.

Denote by $Y^{(1)}\to D\times V$ the bundle of homomorphisms $T_s(D)\to
T_v(V)$. If $D$ and $V$ are given the disk and the almost 
K\"ahler manifold, then
we distinguish the subbundle $X^{(1)}\subset Y^{(1)}$ which consists of
complex linear homomorphisms and we denote $\bar X^{(1)}\to D\times V$ the
quotient bundle $Y^{(1)}/X^{(1)}$. Now, we assign to each $C^1$-map $
v:D\to V$ the section $\bar \partial v$ of the bundle $\bar X^{(1)}$ over
the graph $\Gamma _v\subset D\times V$ by composing the differential of $v$
with the quotient homomorphism $Y^{(1)}\to \bar {X}^{(1)}$. If $c_\delta 
:D\times
V\to \bar X$ is a $H^k-$ section we write $\bar \partial v=c_\delta $ 
for the
equation $\bar \partial v=c_\delta |\Gamma _v$.

\begin{Lemma}
(Gromov\cite{gro})There exists a unique almost complex 
structure $J_g$ on $D\times V$(which 
also depends on the given structures in $D$ and in $V$), such that 
the (germs of) $J_\delta-$holomorphic sections $v:D\to D\times V$ are exactly and 
only the solutions 
of the equations $\bar \partial v=c_\delta $. Furthermore, the 
fibres $z\times V\subset D\times V$ are $J_\delta-$holomorphic(
i.e. the subbundles $T(z\times V)\subset T(D\times V)$ are $J_\delta-$complex) 
and the structure 
$J_\delta|z\times V$ equals the original structure on $V=z\times V$.
\end{Lemma}

\section{Proofs of Theorem1.1}

\begin{Theorem}
There exists a non-constant $J-$holomorphic map $u: (D,\partial D)\to 
(M\times C,W)$ with $E({u})\leq 4\pi r_0^2<s(M,\omega ).$
\end{Theorem}
Proof.  By choose a $J$ such that $s(M,\omega ,J)>4\pi r^2$, we know that 
$J-$holomorphic sphere is removed. 
But the results in 
section 4 shows the solutions of equations (\ref{eq:4.16}) must 
denegerate to a cusp curves, i.e., we obtain a Sacks-Uhlenbeck-Gromov's 
bubble, i.e., disk with boundary 
in $W$ with its area less that $4\pi r_0^2$. 
For more detail, see the proof   
of Theorem 2.3.B in \cite{gro}.

\vskip 3pt 

{\bf Proof of Theorem 1.1}. By Theorem 6.1, 
it is obvious.

\begin{Note}
This work follows from \cite{ma}. One can use the work of this paper to 
reprove or generalizes the results in \cite{ma, ma1}
\end{Note}

\end{document}